\documentclass [11pt,a4paper]{article}
\usepackage{amssymb}
\usepackage{enumerate}
\usepackage{epsfig}
\usepackage{graphicx}
\usepackage{amsmath}
\usepackage{amsfonts}
\usepackage{amsthm}
\usepackage{psfrag}

\newcounter{teoremaganso}

\newtheorem {bigtheo} [teoremaganso] {Theorem}
\newtheorem{thm}{Theorem}[section]
\newtheorem{cor}[thm]{Corollary}
\newtheorem{lem}[thm]{Lemma}
\newtheorem{prop}[thm]{Proposition}
\theoremstyle{definition}

\theoremstyle{remark}
\newtheorem{rem}[thm]{Remark}
\numberwithin{equation}{section}

\newcommand{\R}{I\!\!R}
\newcommand{\N}{I\!\!N}
\newcommand{\eps}{\varepsilon}

\newcommand{\di}{\hbox{div}}
\newcommand{\sgn}{\hbox{sgn}}

\begin{document}
\title{A new uniqueness criterion for the number of periodic orbits of Abel equations}%
\author{M.J. \'{A}lvarez$^\dag$,\,\,\,
  A. Gasull$^+$\,\, and  \,\, H. Giacomini$^\ddag$\\\\
\small $^\dag$\textsl{Departament de Matem\`{a}tiques i
Inform\`{a}tica}
\\
{\small  \textsl{Universitat de les Illes Balears, 07122, Palma de
Mallorca, Spain.}}\\
\small{\texttt chus.alvarez@uib.es}\\
\small $^+$\textsl{Departament de Matem\`{a}tiques, Edifici Cc}
\\
{\small \textsl{Universitat Aut\`{o}noma de Barcelona, 08193
Bellaterra, Barcelona, Spain.}}\\
\small{\texttt gasull@mat.uab.es}\\
\small $^\ddag$\textsl{Laboratoire de Math\'ematiques et Physique
Th\'eorique, CNRS (UMR 6083),}\\
{\small {Facult\'e des Sciences et Techniques}}
\\
{\small \textsl{Universit\'e de Tours, Parc de Grandmont, 37200 Tours, France}}\\
\small{\texttt giacomini@phys.univ-tours.fr}}


\date{}%
\maketitle

\begin{abstract} A solution of the Abel equation $\dot{x}=A(t)x^3+B(t)x^2$ such
that $x(0)=x(1)$ is  called a periodic orbit of the equation. Our
main result proves  that if there exist two real numbers $a$ and
$b$ such that the function $aA(t)+bB(t)$ is not identically zero,
and does not change sign in $[0,1]$ then the Abel differential
equation has at most one non-zero periodic orbit. Furthermore,
when this periodic orbit exists, it is hyperbolic. This result
extends the known criteria about the Abel equation that only refer
to the cases where either $A(t)\not\equiv0$ or $B(t)\not\equiv0$
does not change sign. We apply  this new criterion to study the
number of periodic solutions of two simple cases of Abel
equations: the one where the functions $A(t)$ and $B(t)$ are
1-periodic trigonometric polynomials of degree one and the case
where these two functions are polynomials with three monomials.
Finally, we give an upper bound for the number of isolated
periodic orbits of the general Abel equation
$\dot{x}=A(t)x^3+B(t)x^2+C(t)x,$ when $A(t), B(t)$ and $C(t)$
satisfy adequate conditions.
\end{abstract}

\section{Introduction and main results}\label{section-Main}

In this work we consider smooth Abel equations of the form
\begin{equation}\label{ec-Abel-general}
\dot{x}=\frac{dx}{dt}=A(t)x^3+B(t)x^2,
\end{equation}
defined on the strip  $\mathcal{S}=\{(t,x)\,:\,t\in[0,1]\,\, x\in\mathbb{R}\}.$
 For these equations we study the number of
solutions which start on the line $t=0$ at some point $x=x_0$ and arrive until
the line $t=1$ having the same height, namely $x=x_0.$ These solutions will be
called for short, {\it periodic orbits.} Note that $x=0$ is always a  periodic
orbit of the equation.

Observe that for the Abel equation with  $A$ and $B$ being
1-periodic functions, equation (\ref{ec-Abel-general}) is indeed a
differential equation defined on a cylinder, and this type of {\it
periodic orbits}
 are actual periodic orbits of the Abel equation. In particular,  the trigonometric case
is important because some families of differential systems on the plane can be
transformed, after an adequate change of variables, into this type of Abel
equation, see \cite{Cherkas,Lloyd-JDE}. Thus, the criteria that we obtain in
this paper to control the number of periodic orbits of  Abel equations can also
be used to obtain upper bounds of the number of limit cycles  of several
families of planar polynomial  vector fields.

As usual, we will say that a periodic orbit of equation
(\ref{ec-Abel-general}) is {\it hyperbolic} if the Poincar\'{e}
map between $t=0$ and $t=1$ has derivative different from zero at
the initial condition of the periodic orbit. It is not difficult
to check that $x=0$ is always a non-hyperbolic periodic orbit of
equation (\ref{ec-Abel-general}). We also will say that equation
(\ref{ec-Abel-general}) has a {\it center} at a given periodic
orbit if there exists a neighborhood of this solution where all
the orbits are periodic.

To our knowledge the more general results for bounding the number of periodic
orbits of a complete Abel equation of the form
${dx}/{dt}={A}(t)x^3+{B}(t)x^2+C(t)x,$ are the ones given in
\cite{GL,Lins-Neto,Lloyd} and \cite[Thm. 9.7]{Pliss}. As we will see in
Proposition \ref{nova}, these results applied to the case $C(t)\equiv0$ give
the following result:  If either $A(t)\not\equiv0$ or $B(t)\not\equiv0$ does
not change sign in $[0,1]$ then the maximum number of non-zero periodic orbits
of equation (\ref{ec-Abel-general}) is one. Furthermore, when this periodic
orbit exists it is hyperbolic. Our main theorem extends these results giving a
new criterion of uniqueness of non-zero periodic orbits of equation
(\ref{ec-Abel-general}).

\begin{bigtheo}\label{thm-main}
Consider the Abel equation (\ref{ec-Abel-general}). Assume that
there exist two real numbers $a$ and $b$ such that $aA(t)+bB(t)$
does not vanish identically and does not change sign in $[0,1].$
Then it has at most one non-zero periodic orbit. Furthermore, when
this periodic orbit exists, it is hyperbolic.
\end{bigtheo}

From the proof of the above theorem we can also obtain information about the
location of the non-zero periodic orbit and in some cases prove that this
periodic orbit actually exists, see Remark \ref{existencia}.

In Section \ref{section-Abel-general} we  extend  the above theorem to the more
general Abel  equation
\begin{equation}\label{genbis}
\frac{dx}{dt}={A}(t)x^3+{B}(t)x^2+C(t)x,
\end{equation}
when $\int_0^1 C(t)dt=0.$ We also prove similar, but weaker, results when this
last equality does not hold, see Theorem \ref{thm-con-C}.

It is well known that trigonometrical Abel equations of the form
(\ref{ec-Abel-general}) can have an arbitrary number of isolated
periodic orbits and that this number increases with the degree of
the trigonometrical polynomials $A(t)$ and $B(t);$ see for
instance \cite{Lins-Neto} or \cite{Panov}. For this reason Lins in
\cite{Lins-Neto} and Il'yashenko in \cite{Ilya} have proposed to
study the problem of giving explicit and realistic bounds for
general Abel equations in terms of the degrees of $A$ and $B.$ As
an application of Theorem \ref{thm-main} we consider the simple
case of trigonometrical polynomials of degree one.  The case where
one of the functions has degree zero is much easier and it is
totally solved in Section \ref{pr}, see Remark \ref{constant}.

\begin{bigtheo}\label{thm-trigo-1,1}
Consider the  Abel equation
$$\frac{dx}{dt}=(a_0+a_1 \cos(2\pi t)+a_2\sin(2\pi t))x^3+
(b_0+b_1 \cos(2\pi t)+b_2\sin(2\pi t))x^2,
$$
being $a_0,a_1,a_2,b_0,b_1$ and $b_2$ arbitrary real numbers.
\begin{enumerate}[(a)]
\item It has a center at $x=0$ if and only if
$a_0=b_0=a_2b_1-a_1b_2=0.$

\item If it has not a center at $x=0$ and one of the conditions
$a_0^2\ge a_1^2+a_2^2,$ $b_0^2\ge b_1^2+b_2^2,$ or
$(a_2b_0-a_0b_2)^2+(a_0b_1-a_1b_0)^2\geq (a_2b_1-a_1b_2)^2$ is
satisfied then it has at most one non-zero periodic orbit.
Furthermore, when this periodic orbit exists, it is hyperbolic.

\item There are equations of the above form   having at
least two non-zero hyperbolic periodic orbits.
\end{enumerate}
\end{bigtheo}

On the other hand Abel equations having $A$ and $B$ polynomials in the
$t-$variable are also considered in the literature, see for instance \cite{AL},
\cite{FraYom} or  \cite{Lins-Neto}. We consider the case where $A$ and $B$ are
polynomials with three monomials. From Theorem \ref{thm-main} we obtain the
following result.

\begin{bigtheo}\label{Teorema-Abel-polinomial}
Consider the Abel equation
\begin{equation}\label{Abel-A-B-grados-j-k}
\frac{dx}{dt}=(a_0+a_1t^j+a_2t^k)x^3+(b_0+b_1t^j+b_2t^k)x^2,\end{equation} with
$j,k\in\N$, $0<j<k,$ and $a_0,a_1,a_2,b_0,b_1$ and $b_2$ arbitrary real
numbers.

\begin{enumerate}[(a)]
\item It has a center at $x=0$ if and only if
$a_0+a_1/(j+1)+a_2/(k+1)=b_0+b_1/(j+1)+b_2/(k+1)=a_2b_1-a_1b_2=0.$

\item If it has not a center at $x=0,$ consider the following
sets of conditions:
\begin{enumerate}
\item[1)] either $a_2 b_1-a_1b_2=0 $ or  $\frac{a_2b_0-a_0b_2}{a_2 b_1-a_1b_2}\not\in(-1,0)$,
\item[2)] either $a_1 b_2-a_2b_1=0 $ or $\frac{a_1b_0-a_0b_1}{a_1 b_2-a_2b_1}\not\in(-1,0)$,
\item[3)] either $a_0 b_2-a_2b_0=0 $ or $\frac{a_0b_1-a_1b_0}{a_0 b_2-a_2b_0}\not\in(-1,0)$.
\end{enumerate}

If one of the above conditions  is satisfied, then it has at most
one non-zero periodic orbit. Furthermore, when this periodic orbit
exists, it is hyperbolic.

\item There are equations of the form  (\ref{Abel-A-B-grados-j-k})
having at least two non-zero hyperbolic periodic orbits.

\end{enumerate}

\end{bigtheo}

Notice that the above two theorems do not solve completely  the problem of
bounding the number of periodic orbits for the considered Abel equations. Both
results prove the existence of at most one non-zero hyperbolic periodic orbit,
but under some  hypotheses. On the other hand, in both cases, the maximum
number of non-zero periodic orbits that we have been able to obtain  is two.

Finally, consider the case  where the vector field associated
(\ref{Abel-A-B-grados-j-k}) coincide in both boundaries of $\mathcal{S}$ , {\it
i.e.} when $A(0)=A(1)$ and $B(0)=B(1).$ It is clear that in this case there are
examples of (\ref{Abel-A-B-grados-j-k}) having at least one non-zero periodic
orbit: it suffices to take  for instance $a_1=a_2=b_1=b_2=0$ and $a_0b_0\ne0.$
We have the following corollary which solves the problem of the number of
periodic solutions of this differential equation.

\begin{cor}
Consider equation (\ref{Abel-A-B-grados-j-k}) with the functions
$A$ and $B$ satisfying $A(0)=A(1)$ and $B(0)=B(1),$  {\it i.e.}
$a_1+a_2=b_1+b_2=0$, and not having a center at $x=0$. Then it has
at most one non-zero periodic orbit. Furthermore, when this
periodic orbit exists, it is hyperbolic.
\end{cor}

\section{Preliminary results}\label{pr}
In this section we present some preliminary results. The first one
is   a generalization of the well-known result that gives the
stability of a periodic orbit.  The second one adapts the ideas
presented in \cite{GG} to give upper bounds for the number of
periodic orbits for some planar differential equations to Abel
equations.

Consider a system
\begin{equation}\label{sist-general}
\left\{
\begin{array}{lcl}
\dot{x}&=&P(x,y),\\ \dot{y}&=&Q(x,y)
\end{array} \right.
\end{equation}
of class $N\geq 1,$  and let $L$ be the piece of the orbit associated to a
solution $(x,y)=(\varphi(t),\psi(t))$ for $t\in[0,\tau].$ Assume that through
the two boundary points of $L$ there are defined two transversal curves
$\Sigma_0$ and $\Sigma_\tau.$ By using the continuous dependence with respect
to initial conditions the flow of system (\ref{sist-general}) defines a
Poincar\'{e} map $\Pi$ between $\Sigma_0$ and $\Sigma_\tau.$ To study the
derivative of this map at $L\cap\Sigma_0$ let us parameterize both sections.
Given some small enough $\varepsilon>0,$ consider, without loss of generality,
$\Sigma_t=\{(x,y)=(\alpha(t,n),\beta(t,n))\,:\, |n|<\varepsilon\}$ for $t=0$
and $t=\tau,$ respectively. Furthermore,  assume that
$(\alpha(t,0),\beta(t,0))\in L$ and
$$\Delta(0,0)\Delta(\tau,0)>0,$$
where
\begin{equation}\label{determinante}
\Delta(t,0)=\left|%
\begin{array}{cc}
  \varphi'(t) & \psi'(t)\\
   \frac{\partial \alpha}{\partial n}(t,0) & \frac{\partial \beta}{\partial n}(t,0)
\end{array}%
\right|.
\end{equation}

Following the steps of the study of the stability of a periodic
orbit given in \cite[Sec. 13]{ALGM} we obtain the following
result.

\begin{thm}\label{estabilidad} Let $\Pi$ be the Poincar\'{e} map between two
transversal sections $\Sigma_0$ and $\Sigma_\tau$ of an orbit
$L=\{(\varphi(t),\psi(t))\,:\,t\in[0,\tau]\}$ of
(\ref{sist-general}). Then the derivative of $\Pi$ at $p\in
L\cap\Sigma_0$ is given by
\[
\Pi'(p)=\frac{\Delta(0,0)}{\Delta(\tau,0)}\,\exp\left(\int_0^\tau
\left(\frac{\partial P}{\partial x}(\varphi(s),\psi(s))+\frac{\partial
Q}{\partial y}(\varphi(s),\psi(s))\,\right)ds\right).
\]
where the function $\Delta$ is given in (\ref{determinante}).
\end{thm}

Note that if $L$ is a periodic orbit of (\ref{sist-general}) and
$\tau$ is its period then ${\Delta(0,0)}={\Delta(\tau,0)}$ in the
above formula, giving rise to the well known formula for knowing
the hyperbolicity of a planar periodic orbit.

In the following we will apply the above result to   Abel equations to
determine upper bounds for the number of periodic orbits that they can have on
the strip $\mathcal{S}$.

\begin{cor}\label{corolario}
Let $X(t,x)$ be the  $\mathcal{C}^1$ vector field associated to the system
\begin{equation}\label{ec-Abel}
\dot{t}=1,\quad\quad \dot{x}=A(t)x^3+B(t)x^2+C(t)x, \end{equation} on the strip
$\mathcal{S}$. Let $x=\gamma(t)$ be a solution of (\ref{ec-Abel}) defined in
$[0,1]$. Then, for any non-zero $\mathcal{C}^1$ function $g(t,x)$ the
derivative of the Poincar\'{e} map between $x=0$ and $x=1$ at $(0,\gamma(0))$ is
given by
\[
\Pi'((0,\gamma(0)))=\frac{|g(0,\gamma(0))|}{|g(1,\gamma(1))|}\exp\left(\int_0^{\tau_g}
\mbox{div}(|g(\gamma(t(s)))|X(\gamma(t(s))))\,ds\right),
\]
being $\tau_g$ a positive number and $t(s)$ an increasing function, both
depending on $g$ and given in the proof.
\end{cor}

\begin{proof} Since $g(t,x)$ does not vanish on $\mathcal{S}$  system (\ref{ec-Abel})
is equivalent to
\begin{equation}\label{ec-Abel-g}
\left\{
\begin{array}{llcl}
t'&=dt/ds&=&|g(t,x)|,\\
x'&=dx/ds&=&|g(t,x)|\left(A(t)x^3+B(t)x^2+C(t)x\right).
\end{array}
\right.
\end{equation} Let $(t,\gamma(t))$ be
a solution of (\ref{ec-Abel}).  Then there exists a function
$t(s)$ and a positive number $\tau_g$ such that
$(t(s),\gamma(t(s)))$ is also a solution of (\ref{ec-Abel-g}), for
$s\in[0,\tau_g].$  We take as transversal sections the two borders
of the strip $\mathcal{S}$, parameterized as
$(\alpha(s,n),\beta(s,n))=(t(s),\gamma(t(s))+n)$ for $s=0,\tau_g$.
Then, applying Theorem \ref{estabilidad} we obtain
\[
\Pi'((0,\gamma(0)))=\frac{\Delta(0,\gamma(0))}{\Delta(1,\gamma(1))}\exp\left(\int_0^{\tau_g}
\di(|g(t(s),\gamma(t(s)))|X(t(s),\gamma(t(s))))\,ds\right), \] with
\begin{eqnarray*}
\Delta(t,x)=\left|%
\begin{array}{cc}
  |g(t,x)| & (A(t)x^3+B(t)x^2+C(t)x)|g(t,x)|\\
   0 & 1
\end{array}%
\right|,
\end{eqnarray*}
for $t=0,1$ and then, the result follows.
\end{proof}

Applying the previous corollary  we get the next result.

\begin{cor}\label{Abel}
Consider the Abel equation
\begin{equation}\label{ec-Abel-2}
\frac{dx}{dt}=A(t)x^3+B(t)x^2+C(t)x:=h(t,x), \end{equation} on the strip
$\mathcal{S}$ and let $g(t,x)$ be a non-zero $\mathcal{C}^1$ function,
1-periodic in $t.$ Let $K$ be a connected region where
$\di(|g(t,x)|(1,h(t,x)))$ does not change sign, and vanishes only in a null
measure Lebesgue set, which is not invariant by the flow of equation
(\ref{ec-Abel-2}). Then (\ref{ec-Abel-2}) has, at most, one periodic orbit
completely contained on $K$. Moreover, if it exists, it is hyperbolic and its
stability is given by the sign of $\di(|g(t,x)|(1,h(t,x)))$.
\end{cor}

\begin{proof} Let $x=\gamma(t)$ be a 1-periodic orbit of equation
(\ref{ec-Abel-2}). Then $(t(s),\gamma(t(s)))$ is a solution of
(\ref{ec-Abel-g}), defined for  $s\in[0,\tau_g].$ By using
Corollary \ref{corolario} and the periodicity of $g(t,x)$  we
obtain that the stability of the periodic orbit is given by the
sign of
$$\int_0^{\tau_g}
\di\left(|g(t(s),\gamma(t(s)))|X(t(s),\gamma(t(s)))\right)\,ds. $$

Suppose that there are  several periodic orbits totally contained in the region
$K.$ By using the above result all them have the same stability. Since the
system has no critical points this is not possible and thus it can have at most
one (hyperbolic) periodic orbit fully contained in $K.$

\end{proof}

Now, we can state the result that gives information about the
total number of  periodic orbits of an Abel equation.

\begin{thm}\label{teorema B}
Consider the $\mathcal{C}^1$ Abel equation (\ref{ec-Abel-2}) on the strip
$\mathcal{S}.$ Suppose that there exist $w\in\R$ and a $\mathcal{C}^1$ function
$f(t,x),$ 1-periodic in $t$ and such that
\[
M_w(t,x):=\langle \nabla
f(t,x),\big(1,h(t,x)\big)\rangle+w\,f(t,x)\di\big(1,h(t,x)\big)
\]
does not change sign in the strip, vanishing only in a null
measure Lebesgue set which is not invariant by the flow. Let $K_1$
denote the number of  curves in $\{f=0\}$ joining $t=0$ and $t=1,$
and let $K_2\le K_1$ be the number of these curves which are
invariant by the flow. Then the Abel equation has at most
$K_1+K_2+1$ periodic orbits.

Moreover,  all the orbits that are not contained in $\{f=0\}$ do
not cut this set, are hyperbolic and their stability is given by
the sign of $wf(t,x)M_w(t,x)$ on each of them.
\end{thm}

\begin{proof} Take the value $w$ and the function $f$ given in the statement of
the Theorem.  Consider, for any $\mathcal{C}^1$ function $g,$ the vector field
associated to the Abel equation (\ref{ec-Abel-g}). Note that on $\{f=0\}$ the
function $M_w(t,x)=\langle \nabla f,\big(1,h\big)\rangle$ does not change sign,
which means that the flow of $(1,h(t,x))$ crosses the set $\{f=0\}$ only in one
direction. Then, each periodic orbit either is contained in $\{f=0\}$ or does
not intersect this  set.

In order to bound the number of periodic orbits of the system
which are not contained in $\{f=0\}$  consider in each connected
component of $\mathcal{S}\setminus\{f=0\}$ the function
$g(t,x)=|f(t,x)|^{1/w}$. If
 we compute now $\di(g(t,x)\big(1,h(t,x)\big))$ we get:
\[
\di(g(t,x)\big(1,h(t,x)\big))=\sgn(f)\frac{1}{w}|f|^{\frac{1}{w}-1}M_w(t,x).
\]
In each connected component we can apply Corollary \ref{Abel} and
we get the upper bound stated above.

\end{proof}

\begin{prop}\label{lema-A=cB}
Consider the  differential equation
\[
\frac{dx}{dt}=f(t)P(x),
\]
with $f$ and $P$ smooth functions. Assume that the equation  $P(x)=0$ has
finitely many solutions, $x_1, x_2,\ldots, x_n.$  If $\int_0^1 f(t)\,dt=0$,
then all the solutions $x=x_i$ for $i=1,\ldots,n,$ are centers; otherwise, its
only periodic orbits are  $x=x_i,$  for $i=1,\ldots,n.$ Furthermore, in this
later case, the simple zeros of $P(x)$ are hyperbolic periodic orbits of the
differential equation.
\end{prop}

\begin{proof}

A solution $x=x(t)$ with initial condition $x(0)=\rho\neq x_i,$ for
$i=1,\ldots,n$ satisfies
\begin{equation*}
\varphi(x(t)):=\int_\rho^{x(t)}\frac{du}{P(u)}=\int_0^t f(s)\,ds.
\end{equation*}

Since $\varphi'(x)=\frac{1}{P(x)},$ we have that $\varphi'(x)\neq 0$ if
$\rho\neq x_i$  (because the solutions $x=x_i$ can not be cut). This implies
that $\varphi(x)$ is injective. Then, if $\int_0^1 f(t)\,dt=0$, as
$\varphi(\rho)=0$ and $\varphi(x(1))=0$ we get $x(1)=\rho$ for any solution
$x(t)$ with initial condition $\rho,$ close enough to a periodic orbit $x=x_i.$

On the other hand, if $\int_0^1 f(t)\,dt\neq 0$, we get
$\varphi(x(1))-\varphi(x(0))\neq 0$ and then, $x(1)\neq \rho$ if $\rho\neq x_i$
and the solution $x(t)$ is not periodic.

To prove the hyperbolicity of $x=x_i$ we have to compute the derivative of the
Poincar\'{e} map $\Pi$ defined by the flow between the two transversal sections
$t=0$ and $t=1$ and see that it is not equal to 1. Following \cite{Lloyd} we
get that
\begin{eqnarray*}
\Pi'(x_i)&=&\exp\left(\int_0^1 \frac{\partial}{\partial
x}\Big(P(x)f(t)\Big)\Big|_{x=x_i}\,dt\right)= \exp\left(P'(x_i)\int_0^1
f(t)\,dt\right)\neq 1,
\end{eqnarray*}
as we wanted to prove.
\end{proof}

We also will need some results referred to the case   $C(t)\equiv 0$.

The next lemma is a straightforward consequence of the results of \cite{AL}.

\begin{lem}\label{constants}
Consider the Abel equation (\ref{ec-Abel-general}). The solution $x=0$ is  a
periodic orbit of multiplicity at least two, and the first necessary conditions
to be a center are
\begin{eqnarray*}
V_2&=&\int_0^1B(t)dt=0,\quad V_3=\int_0^1A(t)dt=0 \quad\mbox{ and} \\
V_4&=&\int_0^{1} A(t)\left(\int_0^t B(s)ds \right)dt=0.
\end{eqnarray*}
\end{lem}

\begin{prop}\label{nova}
Consider the Abel equation (\ref{ec-Abel-general}). If either $A(t)\not\equiv0$
or $B(t)\not\equiv0$ does not change sign in $[0,1]$ then the equation has at
most one non-zero periodic orbit and if it exists,  it is hyperbolic.
\end{prop}
\begin{proof} If either $A(t)$ or $B(t)$ does not change sign by using
 the results of \cite{GL} we have that (\ref{ec-Abel-general}) has
 at most three periodic orbits, taking into account their multiplicities.
 By  Lemma \ref{constants} we know that the solution $x=0$
 is at least a double periodic orbit of the Abel equation. Hence we have proved
 that under our hypotheses (\ref{ec-Abel-general}) has at most one non-zero periodic orbit and
  that when it exists, it is hyperbolic.
 \end{proof}

\begin{rem}\label{constant} Note that if we want to study the number of periodic orbits of
equation (\ref{ec-Abel-general}) when either $A(t)$ or $B(t)$ is a constant
function the above two results solve completely the problem, giving the
existence of at most one non-zero hyperbolic orbit. Proposition \ref{lema-A=cB}
solves the case when the constant is zero, while Proposition \ref{nova} solves
the case when it is not zero.
\end{rem}

\section{Lower bounds}\label{lb}

In this section we construct Abel equations with two non-zero hyperbolic
periodic orbits for the two families considered in this paper and we study when
the solution $x=0$ is a center.

\subsection{Trigonometric case}

Consider the  Abel equation
\begin{eqnarray}\label{eq-ejemplo}
\frac{dx}{dt}&=&A(t)x^3+B(t)x^2=(a_0+a_1\cos(2\pi t)+a_2\sin(2\pi t))x^3+\nonumber\\
&&+(b_0+b_1\cos(2\pi t)+b_2\sin(2\pi t))x^2.
\end{eqnarray}
In this section we will see how to  construct examples  with two non-zero
periodic orbits by using two different methods: bifurcating periodic orbits
from $x=0$ and studying perturbations of some centers inside the family.

\subsubsection{Computation of the Lyapunov constants for $x=0$}

The first method we  use to give a lower bound for the number of periodic
orbits of equation (\ref{eq-ejemplo}) consists in computing how many periodic
orbits can bifurcate from  $x=0$. For this purpose we  compute the derivatives
of the Poincar\'{e} map $\Pi$ between $t=0$ and $t=1$, see \cite{AL} or
\cite{Lloyd}. For similarity with the planar case we will call these
derivatives (modulus some non-zero multiplicative constants) the Lyapunov
constants of $x=0.$  We have the following result.
\begin{prop}\label{Lyapunov}
Consider equation (\ref{eq-ejemplo}). The maximum number of
non-zero periodic orbits that can bifurcate from $x=0$ is two.
Moreover, the solution $x=0$ is a center if and only if
$b_0=a_0=a_2b_1-a_1b_2=0.$
\end{prop}

\begin{proof}
We compute the Lyapunov constants for $x=0$ following Lemma
\ref{constants}. The first one is
\[
V_2=\int_0^{1}B(t)\,dt= b_0.
\]
If $V_2\neq 0$, $x=0$ is a semi-stable periodic orbit, {\it i.e.}
if $b_0>0$ the orbits close to $x=0$ with positive initial
condition move away from it and the ones with negative initial
condition approach it, while if $b_0<0$ the orbits with positive
initial condition approach $x=0$ and the ones with negative
initial condition move away.

When $b_0=0$ we compute the next Lyapunov constant:
\[
V_3=\int_0^{1}A(t)\,dt= a_0.
\]
We get similar results as above but now considering $a_0$. If
$a_0=0$, we compute another Lyapunov constant
\[
V_4=\int_0^{1} A(t)\left(\int_0^t B(s)\,ds
\right)\,dt=\frac{a_2b_1-a_1b_2}{4\pi}.
\]
Suppose that this quantity is positive, then if we choose $a_0<0$ and $b_0>0$
with $b_0<<|a_0|<<a_2b_1-a_1b_2$, by a degenerate Hopf bifurcation we can
generate two periodic orbits from $x=0$.

Let us prove now that if $V_2=V_3=V_4=0,$ {\it i.e.} $a_0=b_0=a_2b_1-a_1b_2=0$
then $x=0$ is a center. If $b_1=b_2=0$ we have $B(t)\equiv 0$ and, applying
Proposition \ref{lema-A=cB}  we conclude that $x=0$ is a center. Otherwise,
suppose, for instance, that $b_2\neq 0$. Then $a_1=\frac{a_2b_1}{b_2}$ and we
have $A(t)=\frac{a_2}{b_2}B(t)$. Then, applying again Proposition
\ref{lema-A=cB} with $f(t)=B(t)$ and $P(x)=a_2x^3/b_2+x^2$ we have again that
$x=0$ is a center, because $\int_0^{1} B(t)\,dt=0$.
\end{proof}

\subsubsection{Perturbation of a center}\label{pert-centro-abel}

The second method we use to produce examples with periodic orbits
in equation (\ref{eq-ejemplo}) is the perturbation of a center
inside this family.

\begin{prop}\label{Prop-perturbacion}
Consider the equation
\begin{equation*}
\frac{dx}{dt}= 2\pi \tilde b_1\cos(2\pi t)x^2
+\eps\,\left[(\tilde{a}_0+\tilde{a}_1\cos(2\pi t)+\tilde{a}_2\sin(2\pi t))x^3+
\tilde{b}_0x^2\right],
\end{equation*}
with $\tilde b_1\ne0.$ Then, for $\varepsilon$ small enough, at most two
non-zero periodic orbits bifurcate from the continuous of periodic orbits
existing when $\varepsilon=0.$ Furthermore,  this upper bound can be reached
and the periodic orbits obtained are hyperbolic.
\end{prop}

\begin{proof}  The solutions of the above equation can be
expanded in a small neighborhood of $\eps=0$ as
\[
x_\eps(t;\rho)=x_0(t;\rho)+\eps S(t,\rho)+\eps^2R(t,\rho,\eps),
\]
where $S(0,\rho)=0$ and $x_0(t;\rho)$ is the solution of the
unperturbed equation given by
\begin{equation*}
x_0(t;\rho)=\frac{\rho}{1-\rho\,(\int_0^t B(s)\,ds)}=\frac{
\rho}{1-\rho\,\tilde b_1\sin(2\pi t)}.
\end{equation*}
Set $W(\rho):=S(1,\rho)$. Following the ideas of \cite{Lins-Neto} we know that
the simple zeros of $W(\rho)$ will give rise to initial conditions of periodic
orbits of the perturbed differential equation, which tend to these values when
$\varepsilon$ tends to zero. Doing some computations we get
\begin{eqnarray*}
\hat{W}(\rho):=\frac{W(\rho)}{\rho^2}&=&\int_0^{1}\left(
\tilde{b}_0+\frac{\tilde{a}_0+\tilde{a}_1\cos(2\pi t)+\tilde{a}_2\sin(2\pi
t)}{1-\tilde b_1\rho\sin(2\pi t)}\rho\right)\,dt,
\end{eqnarray*}
which is well defined in the region where the center exists, {\it i.e.}
$|\tilde b_1\rho|<1.$ This fact induces to introduce the natural change of
variable $\tilde b_1\rho=\sin(y),$ for $y\in(-\pi/2,\pi/2),$ for studying the
non-zero zeros of $\hat W.$ We obtain that
\begin{eqnarray*}
\hat{W}(\rho)&=&\tilde b_0+\rho\frac{\tilde a_0\sin(y)+\tilde a_2-\tilde a_2\cos(y)}{\sin(y)\cos(y)}\\
&=& \frac{\tilde a_0\sin(y)+(\tilde b_0\tilde b_1-\tilde a_2)\cos(y)+\tilde
a_2}{\tilde b_1\cos(y)}.
\end{eqnarray*}
Solving the equation $\hat W(\rho)=0$ we get that it has in $(-\pi/2,\pi/2),$
at most two non-zero simple solutions, namely $\rho_1$ and $\rho_2$. These
solutions will give rise, at most, to two non-zero periodic orbits of the
perturbed differential equation, which can be easily seen that are hyperbolic
because $\hat W'(\rho_i)\ne0,$ $i=1,2.$  It is also clear that this upper bound
can be reached. Consider for instance a system with $\tilde a_0=0, \tilde
b_0\tilde b_1-\tilde a_2=1$ and $\tilde a_2=3/4.$ Thus the result follows.
\end{proof}

\subsection{Polynomial case}

In this subsection we consider the  Abel equation (\ref{Abel-A-B-grados-j-k}):

\begin{equation*}
\frac{dx}{dt}=A(t)x^3+B(t)x^2=
        (a_0+a_1t^j+a_2t^k)x^3+(b_0+b_1t^j+b_2t^k)x^2.
\end{equation*}

As in the previous subsection we want to produce examples with as
many non-zero periodic orbits as possible. As before we obtain
examples with two non-zero periodic orbits. However here we only
consider the method of Lyapunov constants, because the other
method gives rise to  tedious integrals and our impression is that
again it does not produce more periodic orbits.

\subsubsection{Computation of the Lyapunov constants for $x=0$}

 We have the
following result.

\begin{prop}\label{cent-pol}
Consider equation (\ref{Abel-A-B-grados-j-k}). The maximum number of non-zero
periodic orbits that can bifurcate from $x=0$ is two. Moreover,  the solution
$x=0$ is a center if and only if
$$a_0+\frac{a_1}{j+1}+\frac{a_2}{k+1}=b_0+\frac{b_1}{j+1}+\frac{b_2}{k+1}=a_2b_1-a_1b_2=0.$$
\end{prop}

\begin{proof} Following the same computations than in Proposition
\ref{Lyapunov} we get that
\[
V_2=\int_0^{1}B(t)\,dt= b_0+\frac{b_1}{j+1}+\frac{b_2}{k+1}.
\]
When $b_0=-((k+1)b_1+(j+1)b_2)/(j+1)(k+1)),$
\[
V_3=\int_0^{1}A(t)\,dt= a_0+\frac{a_1}{j+1}+\frac{a_2}{k+1}.
\]
Finally when $V_2=V_3=0,$
\[
V_4=\int_0^{1} A(t)\left(\int_0^t B(s)\,ds
\right)\,dt=\frac{jk(k-j)(a_2b_1-a_1b_2)}{2(1+j)(2+j)(1+k)(2+k)(2+j+k)}.
\]
If $V_4>0$ we choose $a_0=-\frac{(k+1)a_1+(j+1)a_2}{(j+1)(k+1)}-\mu$ and
$b_0=-\frac{(k+1)b_1+(j+1)b_2}{(j+1)(k+1)}+\lambda$ with $\lambda,\mu>0$, in
such a way that $\lambda<<\mu<<a_2b_1-a_1b_2$. Then, by  a degenerate Hopf
bifurcation we can generate two non-zero periodic orbits from $x=0$. Finally
when $V_2=V_3=V_4=0,$ the fact that the origin is a center follows by
Proposition \ref{lema-A=cB}.
\end{proof}

\section{Proofs of the main results}

\noindent{\it Proof of Theorem \ref{thm-main}}. In our proof we do not care
about the case $ab=0$ because it follows easier than when $ab\neq 0$.
Furthermore, when $ab=0$ the results also follows from \cite{GL}, see
Proposition \ref{nova}.

We start by proving that (\ref{ec-Abel-general}) has at most three non-zero
periodic orbits.

We apply Theorem \ref{teorema B} with $f(x)=x^2(bx-a)$ and $w=-1.$
The function $M_{-1}(t,x)$ is:
\[
M_{-1}(t,x)=x^4(aA(t)+bB(t)),
\]
and it never vanishes identically. Thus by applying Theorem \ref{teorema B} we
know that two types of periodic orbits can exist: the ones contained in the set
$\{f=0\}$, {\it i.e.} $x=a/b$ and $x=0$, and the ones that are contained in
$\mathcal{S}\setminus\{f=0\}$. Since $aA(t)+bB(t)\not\equiv 0$ then $x=a/b$ is
not a periodic orbit.

Thus the  maximum number of non-zero periodic orbits is three, two living in
the same semi-strip where the curve $x=a/b$ is, one bigger than this curve and
the other one smaller, and a third one in the other half-strip. The
hyperbolicity of these periodic orbits is also given by Theorem \ref{teorema
B}.

To prove that in fact there is at most one non-zero periodic orbit
we will compare our differential equation with two simpler ones
for which we know their phase portraits. More precisely we write
the Abel equation (\ref{ec-Abel-general}) in the following two
ways
\begin{eqnarray*}
\frac{dx}{dt}&=&A(t)\,x^2\left(x-\frac{a}{b}\right)+\frac{aA(t)+bB(t)}{b}x^2,\\\\
\frac{dx}{dt}&=&-\frac{b}{a}B(t)\,x^2\left(x-\frac{a}{b}\right)+\frac{aA(t)+bB(t)}{a}x^3.
\end{eqnarray*}

By changing the sign of $t$ if necessary it is not restrictive to assume that
$h(t):=(aA(t)+bB(t))/b\geq 0$. The above equations write as
\begin{eqnarray}
\frac{dx}{dt}&=&A(t)\,x^2\left(x-\frac{a}{b}\right)+h(t)x^2,\label{ec-A}\\\nonumber\\
\frac{dx}{dt}&=&-\frac{b}{a}B(t)\,x^2\left(x-\frac{a}{b}\right)+\frac{b}{a}\,h(t)x^3.\label{ec-B}
\end{eqnarray}

We start comparing (\ref{ec-A}) with
\begin{equation}\label{ec-no-h}
\frac{dx}{dt}=A(t)\,x^2\left(x-\frac{a}{b}\right).
\end{equation}
The global phase portrait of equation (\ref{ec-no-h}) is given by Proposition
\ref{lema-A=cB}. When  $\int_0^1 A(t)dt=0$ the equation has a center at $x=0$
and $x=a/b$ and, since $h(t)x^2\geq 0$, all the solutions of (\ref{ec-no-h}),
except $x=0$, are curves without contact for the flow of (\ref{ec-A}). Thus,
$x=0$ is its only periodic orbit.

If $\int_0^1A(t)dt\neq 0$ then $x=0$ is a double periodic orbit
and $x=a/b$ is a hyperbolic one   and the behavior of the other
solutions is determined by the sign of $-a/b\int_0^1A(t)dt.$
Assume firstly that $-a/b\int_0^1A(t)dt>0$ and that $a/b>0$ (the
case $a/b<0$ follows similarly). Under these inequalities, the
solutions of (\ref{ec-A}) starting below of $x=0$ approach to this
periodic orbit and the solutions starting above $x=0$ approach to
$x=a/b.$ In other words all solutions of (\ref{ec-no-h}),
$x(t;x_0)$ with $x(0;x_0)=x_0$ where $x_0<a/b$ and $x_0\neq 0$
satisfy that $x(1;x_0)>x_0.$ Let $\bar x(t;x_0)$ be the solution
of (\ref{ec-A}) starting also at $x_0.$ Note that  it satisfies
the differential inequality
$$\dot x=A(t)\,x^2\left(x-\frac{a}{b}\right)+h(t)x^2\ge
A(t)\,x^2\left(x-\frac{a}{b}\right),$$ and thus $\bar
x(1;x_0)>x(1;x_0)>x_0$ for all $x_0\ne0$ below $x=a/b.$ Hence,
there is no non-zero periodic orbit below $x=a/b$ and we are done,
because we already know that  in the region $x>a/b$ there is at
most one periodic orbit.

If $-a/b\int_0^1A(t)dt<0$ and  $a/b>0$ the same reasoning can be applied to the
region $x>a/b$ and there is no periodic orbit there. To study the region
$x<a/b$ we distinguish three cases according to the stability of $x=0$ for
(\ref{ec-Abel-general}), which is given by the sign of $\int_0^1B(t)dt$.

If $\int_0^1 B(t)dt=0$, arguing as in the case $\int_0^1A(t)dt=0$
we can prove that $x=0$ is the only periodic orbit. If $\int_0^1
B(t)>0$ there is no periodic orbit in $0<x<a/b$. The reason is
that in this strip there is at most a periodic orbit which, if
exists,  is hyperbolic, and the flow near the boundaries of this
region  implies that if a periodic orbit would exist it should
have even multiplicity.  Thus, the periodic orbit, if exists, it
is unique and  located in the region $x<0$. Finally, if
$\int_0^1B(t)dt<0$ the comparison of the flow of (\ref{ec-A}) and
(\ref{ec-no-h}) and the previous results prove the existence of
exactly one periodic orbit with initial condition between 0 and
$a/b$. The reason is again the sense of the flow in the boundaries
of this region. In order to prove that there is no periodic orbit
in $x<0$, we compare (\ref{ec-A}) with
\begin{equation}\label{ec-no-h-B}
\frac{dx}{dt}=-\frac{b}{a}B(t)\,x^2\left(x-\frac{a}{b}\right).
\end{equation}
Since precisely in the region $x<0$ we have $h(t)x^3\leq 0,$ by using similar
arguments than before we  prove that in this region there is no periodic orbit.
Thus the theorem follows. \hfill $\Box$

\begin{rem}\label{existencia} The proof of Theorem \ref{thm-main} also
helps to locate the non-zero periodic orbit and in some cases to
prove its existence. Collecting all the above results when
$(aA(t)+bB(t))/b\ge0$ and $a/b>0$ we obtain that, if the non-zero
periodic orbit exists, it is located in:
\begin{enumerate}[(i)]
\item The region $x>a/b,$ when $\int_0^1A(t)\,dt<0,$
\item The region $x<0,$ when $\int_0^1A(t)\,dt>0$ and $\int_0^1B(t)\,dt>0,$
\item The region $0<x<a/b$ when $\int_0^1A(t)\,dt>0$ and $\int_0^1B(t)\,dt<0.$   Furthermore in this case
it always exists.
\end{enumerate}
Other signs of $(aA(t)+bB(t))/b$ and $a/b$ can be studied similarly.
\end{rem}

\noindent{\it Proof of Theorem \ref{thm-trigo-1,1}}. (a) Follows
from Proposition \ref{Lyapunov}.

(b) By applying Theorem \ref{thm-main} with $a=-1$ and $b=m$ we get that
\begin{eqnarray*}
M_{-1}(t,x)&=&x^4(mB(t)-A(t))=\\
&=&x^4(b_0 m-a_0+(b_1 m-a_1)\cos(2\pi t)+(b_2 m-a_2)\sin(2\pi t)),
\end{eqnarray*}
and never vanishes identically. Note that if we can find an $m$ such that $(b_0
m-a_0)^2\geq (b_1 m-a_1)^2+(b_2 m-a_2)^2$ then $M_{-1}$ will not change sign.
This $m$ always exists if one of the three conditions of the theorem is
satisfied.  Thus by applying again Theorem \ref{thm-main} we know that there is
at most one non-zero periodic orbit. The hyperbolicity of this periodic orbit
is also given by the same theorem.

(c) The result follows by using either Proposition \ref{Lyapunov}
or Proposition \ref{Prop-perturbacion}.

\hfill $\Box$

\begin{rem}
Notice that  following the proof of Theorem \ref{thm-trigo-1,1} we obtain that
for each $m$ for which $(b_0 m-a_0)^2\geq (b_1 m-a_1)^2+(b_2 m-a_2)^2,$  the
line $x=-1/m$ as a curve without contact for the flow of the differential
equation. Taking all the curves together we obtain a kind of Lyapunov function
that helps to locate the regions where periodic orbits can live. For instance
if  $a_0^2> a_1^2+a_2^2$ or $b_0^2> b_1^2+b_2^2,$ and
$(a_2b_0-a_0b_2)^2+(a_0b_1-a_1b_0)^2< (a_2b_1-a_1b_2)^2$ then all the non zero
values of $m$ are allowed in the proof of the above theorem and we can show
that there is no non-zero periodic orbit.
\end{rem}

\noindent{\it Proof of Theorem \ref{Teorema-Abel-polinomial}.} (a)
Follows from Proposition \ref{cent-pol}.

(b) The proof follows in a very similar way than the proof of (b) of Theorem
\ref{thm-trigo-1,1}. In order to apply Theorem \ref{thm-main} we set $a=-b_2$
and $b=a_2$. Thus,
\[
M_{-1}(t,x)=x^4(a_2B(t)-b_2A(t))=x^4(a_2b_0-a_0b_2+(a_2b_1-a_1b_2)t^j),
\]
and, if condition {\it 1)} is satisfied, $M_{-1}(t,x)$ has a definite sign in
all the strip.

To study the differential equation under condition {\it 2)} we can
consider $a=-b_1$ and $b=a_1$. Then the function $M_{-1}(t,x)$ is
given by
\[
M_{-1}(t,x)=x^4(a_1b_0-a_0b_1+(a_1b_2-a_2b_1)t^k).
\]

Finally, when condition {\it 3)} is considered, by using $a=-b_0$ and $b=a_0$,
we get that
\[
M_{-1}(t,x)=x^4t^j(a_0b_1-a_1b_0+(a_0b_2-a_2b_0)t^{k-j}).
\]

(c) The result follows by using  Proposition \ref{cent-pol}.

\hfill $\Box$

\section{The general Abel equation}\label{section-Abel-general}

Regarding to the general Abel equation (\ref{genbis}):
\begin{equation*}
\frac{dx}{dt}={A}(t)x^3+{B}(t)x^2+C(t)x,
\end{equation*}
we have proved the following two extensions of Theorem \ref{thm-main}:

\begin{thm} Consider the general Abel equation (\ref{genbis}).
Assume that $\int_0^1 C(t)\,dt=0$ and that there exist two real numbers $a$ and
$b$ such that
$$aA(t)\exp\left(\int_0^t C(s)\,ds\right)+bB(t)$$ does not vanish identically
and does not change sign  for all $t\in[0,1]$. Then it has at most
one non-zero periodic orbit. Furthermore, when this periodic orbit
exists, it is hyperbolic.
\end{thm}

\begin{proof} If $\int_0^1 C(t)dt=0$, the well known change of variables $$y=x\exp\left(-\int_0^t
C(s)\,ds\right)$$ transforms equation  (\ref{genbis}) into
\begin{equation*}
\frac{dy}{dt}={A}(t)\exp\left(2\int_0^t
C(s)\,ds\right)y^3+{B}(t)\exp\left(\int_0^t C(s)\,ds\right)y^2.
\end{equation*}
Furthermore this change sends periodic orbits of (\ref{genbis}) into periodic
orbits of the above equation. By applying Theorem \ref{thm-main} to this
equation the result follows.
\end{proof}

\begin{thm}\label{thm-con-C}
Consider equation (\ref{genbis}). Assume that there exist three real numbers
$a,b$ and $c$ such that $a {A}(t)+b {B}(t)$ is not identically zero and does
not change sign for all $t\in[0,1],$ and that $(bC(t)-c {A}(t))^2+(a {A}(t)+b
{B}(t))(c {B}(t)+aC(t))<0$ for all $t\in[0,1].$ Then it has at most four
non-zero periodic orbits.
\end{thm}

\begin{proof}
As in the proof of Theorem \ref{thm-main} we can only consider the case
$ab\ne0$ and we apply Theorem \ref{teorema B} with $f(x)=bx^3-ax^2+cx$ and
$w=-1.$ The function $M_{-1}(t,x)$ is:
\[
M_{-1}(t,x)=x^2\Big((a {A}(t)+b {B}(t))x^2+2(bC(t)-c {A}(t))x-(c
{B}(t)+aC(t))\Big).
\]
If the two conditions of the statement of the theorem are satisfied then this
function does not change sign, and by using Theorem \ref{teorema B}, we know
that only two types of periodic orbits can exist: the ones contained in the set
$\{f=0\}$, {\it i.e.} $x^\pm:=(a\pm\sqrt{a^2-4bc})/(2b)$ and $x=0$, and the
ones that are contained in $\mathcal{S}\setminus\{f=0\}$. Furthermore, there is
at most one periodic orbit of this type in each connected component of
$\mathcal{S}\setminus\{f=0\}$.

Thus,  when $\Delta:=a^2-4bc<0$ ($\Delta=0,\Delta>0,$ respectively ) equation
(\ref{genbis}) has at most 2 (4, 6, respectively) non-zero periodic orbits.
Hence, when $a^2-4bc\le0$ we are done. When $a^2-4bc>0$ and none of the lines
$x=x^{\pm}$ is a periodic orbit  the result also follows. Finally assume that one
of these lines is a periodic orbit. In this situation  it is easy to see that
the other one can not be  a periodic orbit of the equation. Then we have proved
that the maximum number of non-zero periodic orbits is five: the one given by
the invariant line, say $x=x^+\ne0,$ and the other four contained in
$\mathcal{S}\setminus\{f=0\}$. To end the proof we will see that in this case
there are at most three periodic orbits in $\mathcal{S}\setminus\{f=0\}$. If
$x=x^+$ is a periodic orbit of (\ref{genbis}) then
\[
0\equiv\frac{dx}{dt}\Big|_{x=x^+}=x^+\,\left( {A}(t) x^+\,(x^+ -\frac{a}{b})
+C(t) + x^+h(t)\right),
\]
where $h(t)=(a {A}(t)+b {B}(t))/b$. Then $C(t)=
{A}(t)x^+(\frac{a}{b}-x^+)-x^+h(t)$ and the equation (\ref{genbis}) writes as
\[
\frac{dx}{dt}=x\,(x -x^+) \,\left( {A}(t)(x-x^-)+h(t)\right).
\]

We will compare  the solutions of the previous equation with the solutions of
\[
\frac{dx}{dt}=A(t)x\,(x -x^+)\left(x-x^-\right).
\]
By using  Proposition \ref{lema-A=cB} we know that when $\int_0^1 A(t)\,dt=0$
this equation has a center at $x=0,x=x^+$ and $x=x^-,$ respectively. When
$\int_0^1 A(t)\,dt\ne0,$ $x=0,$ $x=x^+$ and $x=x^-$ are the only  periodic
orbits. Following the same reasoning than in the proof of Theorem
\ref{thm-main}, we can distinguish the cases according to the sign of several
involved functions, obtaining in all cases  that there is no periodic orbit in
one of the four connected components of the set $\mathcal{S}\setminus\{f=0\}.$
Thus the result follows.
\end{proof}

\bibliographystyle{amsplain}

\end{document}